\documentclass[11pt]{article}
\usepackage{amsfonts}
\usepackage{latexsym,amsmath}
\usepackage{amssymb,array}
\usepackage[sort&compress,round,comma]{natbib}
\usepackage{setspace}
\usepackage{enumitem}
\usepackage{graphicx}
\usepackage{xcolor}
\usepackage{float}
\setlist{noitemsep}
\makeatletter \oddsidemargin 0in \evensidemargin 0in \textwidth 16cm
\usepackage{graphicx}
\providecommand{\keywords}[1] {
  \textit{Keywords:} #1}
\begin{document}
    \newcommand{\ov}{\overline}
    \newcommand{\om}{\omega}
    \newcommand{\ga}{\gamma}
    \newcommand{\cd}{\circledast}
    \newtheorem{thm}{Theorem}[section]
    \newtheorem{remark}{Remark}[section]
    \newtheorem{counterexample}{Counterexample}[section]
    \newtheorem{coro}{Corollary}[section]
    \newtheorem{propo}{Proposition}[section]
    \newtheorem{definition}{Definition}[section]
    \newtheorem{example}{Example}[section]
    \newtheorem{lem}{Lemma}[section]
    \numberwithin{equation}{section}
    \date{}

\title{Some stochastic comparison results for frailty and resilience models}
\author{Arindam Panja$^1$, Pradip Kundu$^2$\footnote
    {Corresponding author}
     and Biswabrata Pradhan$^1$\\
$^1$ SQC \& OR Unit, Indian Statistical Institute, Kolkata-700108,
India\\ $^2$ Decision Sciences, Birla Global University,
Bhubaneswar-751003, India\let\thefootnote\relax\footnote{Email
Address: arindampnj@gmail.com (Arindam Panja), kundu.maths@gmail.com
(Pradip Kundu), bis@isical.ac.in (Biswabrata Pradhan)}} \maketitle
\begin{abstract}
Frailty and resilience models provide a way to introduce random
effects in hazard and reversed hazard rate modeling by random
variables, called frailty and resilience random variables,
respectively, to account for unobserved or unexplained heterogeneity
among experimental units. This paper investigates the effects of
frailty and resilience random variables on the baseline random
variables using some shifted stochastic orders based on some ageing
properties of the baseline random variables. Relevant examples are
provided to illustrate the results. Some results are illustrated
with real-world data.
\end{abstract}

\section{Introduction:}
Heterogeneity is a very common issue in many areas including
reliability, survival analysis, demography and epidemiology. For
instance, in mechanical systems, heterogeneity occurs due to
unit-to-unit variability, changes in operating environments, the
diversity of tasks and workloads during its lifetime. For example,
as mentioned in \cite{ozek} that a complex device like an airplane
has large number of components where the failure structure of each
component depends on a set of environmental conditions (e.g. the
levels of vibration, atmospheric pressure, temperature, etc.) that
vary during take-off, cruising and landing. So incorporating
heterogeneity into hazard (failure) rate modeling is a common
practice to achieve accuracy in the estimation. The proportional
hazard (PHR) model is the most applied model in the case where
factors (covariates) influencing the environment/operating condition
are known and can be quantified. In such case, hazard rate of an
individual is considered to be constant multiplicative to the
baseline hazard. However, in many practical situation it may happen
that some factors influencing the operating condition are unknown,
and heterogeneity occurs in an unpredicted and unexplained manner. A
component may subject to different levels of operating environment
(e.g. voltage, temperature) which is not fixed but changes over
time. Component lifetimes and reliability depend on these random
environmental variations. Frailty models
(\citealp{cha,da,houg,Gupt2,guptp,li,vaup,zaki}) provide a way to
introduce random effects in the model by a random variable (r.v.),
called frailty r.v., to account for unobserved (unexplained)
heterogeneity among experimental units in their hazard (failure)
rates. For instance, \cite{vaup} discussed that in survival
analysis, mortality of individuals differ due to large number of
factors beyond age, e.g. the individual's susceptibility to causes
of death, response to treatment and various risk factors. They
considered a frailty r.v. to cope with the unobserved individual
differences in mortality rates while defining the force of mortality
of individuals. \cite{chah} considered a frailty r.v. in the model
for mission abort/continuation policy for heterogeneous systems, to
justify heterogeneity which may occur due to various reasons such as
quality of resources used in the production process, operation and
maintenance history, and human errors.
\par Let $X$ be a r.v. with distribution function $F$ and survival function $\bar{F}=1-F$, and $\Lambda$ be a continuous r.v. with distribution function $H$, probability density
function (pdf) $h$ and hazard function $r_{F}$. A r.v. $X^{*}$ is
said to follow multiplicative frailty model with baseline
distribution $F$ and frailty r.v. $\Lambda$ if its survival function
is given by
\begin{equation}\label{frlt}
    \bar{F}^*(t)= \int_{0}^{\infty}\bar{F}^{\lambda}(t)  dH(\lambda)
\end{equation}
Here the frailty r.v. $\Lambda$ serves as an unobserved random
factor that modifies multiplicatively the underlying hazard function
$r_{F}$ of an individual such that the individual is supposed to
have hazard rate $r_{F^*}(t)$ at age $t$, so that given
$\Lambda=\lambda$, the conditional hazard rate function of $X^{*}$
will be $r_{F^*}(t|\lambda)=\lambda~ r_{F}(t)$, $t\geq 0$.

\par In analogy to the frailty model, to account for unobserved/unexpalined heterogeneity in
the reversed hazard rates of the experimental units, the resilience
model (reversed frailty models) is introduced. A r.v. $X^{*}$ is
said to follow resilience model with baseline distribution $F$ and
resilience r.v. $\Lambda$ if its distribution function is given by
\begin{equation}\label{rslng}
    F^*(t)= \int_{0}^{\infty}F^{\lambda}(t)  dH(\lambda)
\end{equation}

The model (\ref{frlt}) is also regarded as mixture (continuous)
distribution of the PHR model with baseline distribution function F,
and mixing r.v. $\Lambda$ (\citealp{da}). Similarly model
(\ref{rslng}) is regarded as the mixture distribution of the
proportional reversed hazard model (\citealp{li}).

Ageing properties and stochastic comparisons of frailty models,
arising from different choices of frailty/baseline distributions,
have been studied by \cite{Gupt1,kaya,kayb,misra,xie} and \cite{Xu}.
On the other hand ageing properties and stochastic comparison of
resilience models have been studied by \cite{Gupta2} and \cite{li}
considering different baseline distributions and/or resilience
distributions. \cite{he} derived comparison results for general
weighted frailty models with respect to some relative stochastic
orders. Another important study in this area is to compare $X^{*}$
and $X$ which helps us to understand the effect of
frailty/resilience r.v. on the underlying original (baseline)
distribution. \cite{misra} studied the ageing of r.v. following
frailty (resilience) model ($X^{*}$) relative to the ageing of
corresponding baseline r.v. $X$. \cite{da} compared $X^{*}$ with
certain frailty, i.e. $X^{*}|\Lambda=\lambda$ and $X$ satisfying
some conditions on the mean of the frailty.
\par In our work we study the effects of frailty and resilience r.v.'s on the baseline r.v. ($X$) using some shifted stochastic orders based on some ageing
properties of $X$. In Section 2, we provide definitions of sifted
stochastic orders along with their usefulness in stochastic
comparisons and superiority over their usual versions. In Section 3,
we study the effect of frailty r.v. on the baseline r.v., i.e.
compare $X^{*}$ and $X$ with respect to some shifted stochastic
orders, where in Section 4, a similar study is carried out in case
of resilience model. In Section 5, we illustrate some of our derived
results with real-world data.
\section{Notations and definitions}
Let $X$ and $Y$ be two absolutely continuous and non-negative r.v.'s
with distribution functions $F$ and $G$; survival (/reliability)
functions $\bar{F}$ and $\bar{G};$ pdf $f$ and $g$; hazard (failure)
rate functions $r_X$ and $r_Y$; reversed hazard rate functions
$\tilde{r}_X$ and $\tilde{r}_Y$, respectively.

\begin{definition}

    $X$ is said to be smaller than $Y$ in the
    \begin{enumerate}
        \item [(a)] up (down) shifted hazard rate order denoted by $X\leq_{hr\uparrow} Y$ ($X\leq_{hr\downarrow} Y$), if $\bar{G}(x) \Big/ \bar{F}(x+t)$ $\left(\bar{G}(x+t) \Big/ \bar{F}(x)\right)$ is increasing in $x$ for all
        $t>0$ (\citealp{lillo}).\vspace{0.08cm}
        \item [(b)] up (down) shifted reversed hazard rate order denoted by $X\leq_{rh\uparrow} Y$ ($X\leq_{rh\downarrow} Y$), if $G(x) \Big/ F(x+t)$ $\left(G(x+t) \Big/ F(x)\right)$ is increasing in $x$ for all
        $t>0$ (\citealp{di}).\vspace{0.08cm}
        \item [(c)] up (down) shifted likelihood ratio order denoted by $X\leq_{lr\uparrow} Y$ ($X\leq_{lr\downarrow} Y$), if $g(x) \Big/ f(x+t)$ $\left(g(x+t) \Big/ f(x)\right)$ is increasing in $x$ for all $t>0$ (\citealp{lillo,shaked}).\vspace{0.08cm}
        \item [(d)] up (down) shifted mean residual life order denoted by $X\leq_{mrl\uparrow} Y$ ($X\leq_{mrl\downarrow} Y$), if $\int_{x+t}^{\infty}\bar{G}(u)du \Big/\int_{x}^{\infty} \bar{F}(u)du$ $\left(\int_{x}^{\infty}\bar{G}(u)du \Big/\int_{x+t}^{\infty} \bar{F}(u)du\right)$ is increasing in $x$ for all $t>0$ (\citealp{nanda}).\vspace{0.08cm}
        \item [(e)] up (down) shifted mean inactivity time order (also known as reversed mean residual life order) denoted by $X\leq_{mit\uparrow} Y$ ($X\leq_{mit\downarrow} Y$), if $\int_{0}^{x+t}F(u) \Big/ \int_{0}^x G(u) du$ $\left(\int_{0}^{x}F(u) \Big/ \int_{0}^{x+t} G(u) du \right)$ is decreasing in $x$ for all $t>0$ (\citealp{nandab,kayb}).
    \end{enumerate}
\end{definition}
It is worth to mention that $X\leq_{lr\uparrow} Y\Rightarrow
X\leq_{lr} Y$, $X\leq_{hr\uparrow} Y\Rightarrow X\leq_{hr} Y$,
$X\leq_{rh\uparrow} Y \Rightarrow X\leq_{rh} Y$, and
$X\leq_{mrl\uparrow} Y \Rightarrow X\leq_{mrl} Y$. For positive
support, these results also hold true for respective down shifted
orders. That means these shifted stochastic orders are stronger than
their respective usual versions of stochastic orders. Also, these
shifted orders can be considered as generalization of their usual
counterparts in some aspects. For instance, unlike likelihood ratio
ordering, shifted likelihood ratio ordering preserves the order
under convolution (\citealp{lillo}). If $X\leq_{lr\uparrow} Y$, then
$\kappa_X (t_1)\leq \kappa_Y (t_2)$ for $t_1 \geq t_2 \geq 0$, where
$\kappa_X\equiv f'/f$ and $\kappa_Y\equiv g'/g$ (\citealp{lillo}).
Note that if $X\leq_{lr} Y$, then $\kappa_X (t)\leq \kappa_Y (t)$
for all $t\geq 0$. It is shown by \cite{di} and \cite{lillo} that
$X\leq_{hr\uparrow} Y \Longleftrightarrow r_X (t_1)\geq r_Y (t_2)$
for $t_1 \geq t_2\geq 0$. Note that $X\leq_{hr} Y$ implies $r_X
(t)\geq r_Y (t)$ for all $t\geq 0$. Similarly, $X\leq_{rh\uparrow} Y
\Longleftrightarrow \tilde{r}_X (t_1)\leq \tilde{r}_Y (t_2)$ for
$t_1 \geq t_2\geq 0$ (\citealp{di}). Note that $X\leq_{rh} Y$
implies $\tilde{r}_X (t)\leq \tilde{r}_Y (t)$ for all $t\geq 0$.
Also if $X\leq_{rh\uparrow} Y$, then $\bar{F}(t_1)\leq \bar{G}(t_2)$
for $t_1\geq t_2 \geq 0$. If $X\leq_{mrl\uparrow} Y$, then
$m_X(t_1)\leq m_Y(t_2)$ for $t_1\geq t_2\geq 0$, where
$m_X(t)=\int_{t}^{\infty}\bar{F}(u)du \Big/\bar{F}(t)$ is the mean
residual life (mrl) of $X$ (\citealp{nanda}). If
$X\leq_{mit\uparrow} Y$, $mit_{X} (t_1)\geq mit_Y (t_2)$ for $t_1
\geq t_2\geq 0$, where $mit_X(t)=\int_{0}^t F(u) du\Big/F(t)$ is
known as mean inactivity time (or reversed mean residual life) of
$X$. Similar results are also shown for down shifted orders, e.g.,
if $X \leq_{hr\downarrow} Y$, then $r_X (t_1)\geq r_Y (t_2)$ for
$t_2 \geq t_1\geq 0$ (\citealp{lillo}). Thus these shifted
stochastic orders give us the flexibility that even at different
points of time for the two variables, we can compare their hazard
rates, reversed hazard rates, survival functions, mean residual
life, etc. One such specific instance is that we can compare the
reliability of an used device and a new device using the shifted
stochastic orders. For more discussion on those shifted orders
including their applications and preservations properties, we refer
to \cite{abo,naq,kayb} and references therein.
\par Next we give the definitions of some ageing classes
(\citealp{laim}).
\begin{definition}
    $X$ is said to have
    \begin{enumerate}
        \item [(a)] increasing (decreasing) likelihood ratio (ILR (DLR)) if $f$ is log-concave (log-convex)  or equivalently for any $t>0$, $f(x+t)/f(x)$ is decreasing (increasing) in $x.$
        \item[(b)] increasing (decreasing) failure rate (IFR (DFR)) if $\bar{F}$ is log-concave (log-convex) or equivalently for any $t>0$, $\bar{F}(x+t)/\bar{F}(x)$ is decreasing (increasing) in $x.$
        \item[(c)] decreasing (increasing) reversed failure rate (DRFR (IRFR)) if $F$ is log-concave (log-convex)  or equivalently for any $t>0$, $F(x+t)/F(x)$ is decreasing (increasing) in $x.$
        \item [(d)] decreasing (increasing) mean residual life (IMRL(DMRL)) if $\int_{x}^{\infty} \bar{F}(u) du$ is log-convex (log-concave) or equivalently for any $t>0$, $\int_{x+t}^{\infty} \bar{F}(u) du / \int_{x}^\infty \bar{F}(u)du$ increasing (decreasing) in $x.$
        \item[(e)] increasing mean inactivity time (IMIT) if $\int_{0}^{x} F(u) du$ is log-concave or equivalently for any $t>0$, $\int_{0}^{x+t} F(u) du / \int_{0}^x F(u)du$  decreasing in $x.$
    \end{enumerate}
\end{definition}

\section{Results for frailty model}
Here we study effect of frailty r.v. on the baseline r.v. with
respect to some shifted stochastic ordering based on some ageing
properties of concerned baseline r.v.'s. Throughout this section, we
consider $X$ and $X^*$ be two r.v.'s as defined in Section 1 for
which the survival function of $X^*$ is given by the equation
(\ref{frlt}). Also consider that $X$ be an absolutely continuous
non-negative r.v.
\par In the following theorem we derived that, for a baseline r.v.
$X$ with ILR (resp. DLR) property, effect of a frailty r.v.
$\Lambda$ with $P(0<\Lambda\leq 1)=1$ (resp. $P(\Lambda\geq 1)=1$)
on $X$ is that, $X^*$ will be greater than (resp. less than) $X$ in
the sense of the up shifted likelihood ratio order.
\begin{thm}\label{th1}~
    \begin{enumerate}
        \item[(i)] $X^* \geq_{lr\uparrow} X$ if $X$ is ILR, provided $0<\Lambda\leq 1$ with probability 1;
        \item[(ii)] $X^* \leq_{lr\uparrow} X$ if $X$ is DLR, provided $\Lambda\geq 1$ with probability 1.
    \end{enumerate}
\end{thm}
\textbf{Proof:} \begin{enumerate}

\item[(i)] We have \begin{eqnarray}
\nonumber \frac{f^*(x)}{f(x+t)} &=& \frac{f(x)}{f(x+t)}\times
\int_{0}^{\infty} \lambda \bar{F}^{\lambda-1}(x) dH(\lambda)
\\&=&\label{ulr1} \displaystyle E\left[\frac{f(x)\Lambda \bar{F}^{\Lambda-1}(x)}{f(x+t)}\right]
\end{eqnarray}

Now $X$ is ILR implies $\frac{f(x)}{f(x+t)}$ is increasing in $x$
for any $t>0$. Again  $\lambda\bar{F}^{\lambda -1}(x)$ will be
increasing in $x$ for any $0<\lambda\leq 1.$ Now if we consider
$\Lambda$ such that $P(0<\Lambda\leq 1)=1$ the result follows
immediately.
\item[(ii)] Similarly $X$ is DLR implies $\frac{f(x)}{f(x+t)}$ is decreasing in $x$. Again  $\lambda\bar{F}^{\lambda -1}(x)$ will be decreasing in $x$ for any $\lambda\geq 1.$ Now if we consider $\Lambda$ such that $P(\Lambda\geq1)=1$, the result follows immediately.
\end{enumerate}

Examples \ref{ex1} and \ref{ex2} illustrate Theorem \ref{th1}(i) and
\ref{th1}(ii) respectively.
\begin{example}\label{ex1}
Let $X$ be a  gamma r.v. with pdf $f(x)=xe^{-x},~x\geq 0.$ Then
clearly $X$ is ILR. Consider the frailty r.v. $\Lambda$ to be
uniformly distributed on $[0,1].$ Then it is easy to check that
$f^*(x)/f(x+t)$ is increasing in $x$ for all $t>0$, giving $X
\leq_{lr\uparrow} X^*$.
\end{example}

\begin{example}\label{ex2}
Let $X$ be a  Weibull r.v. with pdf $f(x)=3x^2e^{-x^3},~x\geq 0.$
Then clearly $X$ is ILR. Consider the frailty r.v. $\Lambda$ to be
uniformly distributed on $[1,3].$ Then it is easy to check that
$f^*(x)/f(x+t)$ is decreasing in $x$ for all $t>0$, giving $X
\geq_{lr\uparrow} X^*$.
\end{example}

The following corollary follows immediately in case $\Lambda$ is a
degenerate r.v.
\begin{coro}~
    \begin{enumerate}
        \item[(i)] $X^* \geq_{lr\uparrow} X$ if $X$ is ILR, provided $0<\lambda\leq 1$ ;
        \item[(ii)] $X^* \leq_{lr\uparrow} X$ if $X$ is DLR, provided $\lambda\geq 1$ .
    \end{enumerate}
\end{coro}

\begin{thm} \label{th2}~
    \begin{enumerate}
        \item[(i)] $X^* \geq_{lr\downarrow} X$ if $X$ is DLR, provided $0<\Lambda\leq 1$ with probability 1;
        \item[(ii)] $X^* \leq_{lr\downarrow} X$ if $X$ is ILR, provided $\Lambda\geq 1$ with probability 1.
    \end{enumerate}
\end{thm}
\textbf{Proof:} \begin{enumerate}

    \item[(i)] We have \begin{eqnarray}
        \nonumber \frac{f^*(x+t)}{f(x)} &=& \frac{f(x+t)}{f(x)}\times \int_{0}^{\infty} \lambda \bar{F}^{\lambda-1}(x+t) dH(\lambda)
        \\&=&\label{dlr1} \displaystyle E\left[\frac{f(x+t)\Lambda \bar{F}^{\Lambda-1}(x+t)}{f(x)}\right]
    \end{eqnarray}

    Now $X$ is DLR implies $\frac{f(x+t)}{f(x)}$ is increasing in $x$ for any $t>0$. Again  $\lambda\bar{F}^{\lambda -1}(x+t)$ will be increasing in $x$ for any $0<\lambda\leq 1.$ Now if we consider $\Lambda$ such that $P(0<\Lambda\leq 1)=1$ the result follows immediately.
    \item[(ii)] Similarly $X$ is ILR implies $\frac{f(x+t)}{f(x)}$ is decreasing in $x$. Again  $\lambda\bar{F}^{\lambda -1}(x+t)$ will be decreasing in $x$ for any $\lambda\geq 1.$ Now if we consider $\Lambda$ such that $P(\Lambda\geq 1)=1$ the result follows immediately.
\end{enumerate}
\begin{remark} Theorem \ref{th2}(i) implies that under the stated
assumptions on $X$ and $\Lambda$, $\kappa_{X^*}(t)\geq \kappa_X
(t')$ for $t\geq t' \geq 0$. Similarly, Theorem \ref{th2}(ii)
implies that $\kappa_{X^*}(t)\leq \kappa_X (t')$ for $t'\geq t \geq
0$.\end{remark}

The following theorem shows that, for a baseline r.v. $X$ with IFR
(resp. DFR) property, effect of a frailty r.v. $\Lambda$ with
$P(0<\Lambda\leq 1)=1$ (resp. $P(\Lambda\geq 1)=1$) on $X$ is that,
$X^*$ will be greater than (resp. less than) $X$ in the sense of the
up shifted hazard rate order.
\begin{thm}\label{th3}~
    \begin{enumerate}
        \item[(i)] $X^* \geq_{hr\uparrow} X$ if $X$ is IFR, provided $0<\Lambda\leq 1$ with probability 1;
        \item[(ii)] $X^* \leq_{hr\uparrow} X$ if $X$ is DFR, provided $\Lambda\geq 1$ with probability 1.
    \end{enumerate}
\end{thm}
\textbf{Proof:} \begin{enumerate}

    \item[(i)] We have \begin{eqnarray}
        \nonumber \frac{\bar{F}^*(x)}{\bar{F}(x+t)} &=& \frac{\bar{F}(x)}{\bar{F}(x+t)}\times \int_{0}^{\infty}  \bar{F}^{\lambda-1}(x) dH(\lambda)
        \\&=&\label{dhr2} \displaystyle E\left[\frac{\bar{F}(x) \bar{F}^{\Lambda-1}(x)}{\bar{F}(x+t)}\right]
    \end{eqnarray}

    Now $X$ is IFR implies $\frac{\bar{F}(x)}{\bar{F}(x+t)}$ is increasing in $x$ for any $t>0$. Again  $\displaystyle\bar{F}^{\lambda -1}(x)$ will be increasing in $x$ for any $0<\lambda\leq 1.$ Now if we consider $\Lambda$ such that $P(0<\Lambda\leq 1)=1$ the result follows immediately.
    \item[(ii)] Similarly $X$ is DFR implies $\frac{\bar{F}(x)}{\bar{F}(x+t)}$ is decreasing in $x$. Again  $\bar{F}^{\lambda -1}(x)$ will be decreasing in $x$ for any $\lambda\geq 1.$ Now if we consider $\Lambda$ such that $P(\Lambda\geq 1)=1$ the result follows immediately.
\end{enumerate}
\begin{remark} Theorem \ref{th3}(i) implies that under the stated
assumptions on $X$ and $\Lambda$, $r_{X^*}(t)\leq r_X (t')$ for
$t'\geq t \geq 0$. Similarly, Theorem \ref{th3}(ii) implies that
$r_{X^*}(t)\geq r_X (t')$ for $t\geq t' \geq 0$.\end{remark}

Examples \ref{ex3} and \ref{ex32} illustrate Theorem \ref{th3}(i)
and \ref{th3}(ii) respectively.
\begin{example}\label{ex3}
Let $X$ follows Weibull distribution with sf $\bar{F}(x)= e^{-x^2}$,
$x\geq0$. Clearly, $X$ is IFR. Let the frailty r.v. $\Lambda$ to be
uniformly distributed on $[0,1].$ Then it is easy to check that
$\bar{F}^*(x)/\bar{F}(x+t)$ is increasing in $x$ for all $t>0.$
\end{example}

\begin{example}\label{ex32}
Let $X$ follows Weibull distribution with sf $\bar{F}(x)=
e^{-x^{0.5}}$, $x\geq0$. Clearly, $X$ is DFR. Let the frailty r.v.
$\Lambda$ to be uniformly distributed on $[2,5].$ Then it is easy to
check that $\bar{F}^*(x)/\bar{F}(x+t)$ is decreasing in $x$ for all
$t>0.$
\end{example}

\begin{thm}\label{th4}~
    \begin{enumerate}
        \item[(i)] $X^* \geq_{hr\downarrow} X$ if $X$ is DFR, provided $0<\Lambda\leq 1$ with probability 1;
        \item[(ii)] $X^* \leq_{hr\downarrow} X$ if $X$ is IFR, provided $\Lambda\geq 1$ with probability 1.
    \end{enumerate}
\end{thm}
\textbf{Proof:} \begin{enumerate}

    \item[(i)] We have \begin{eqnarray}
        \nonumber \frac{\bar{F}^*(x+t)}{\bar{F}(x)} &=& \frac{\bar{F}(x+t)}{\bar{F}(x)}\times \int_{0}^{\infty}  \bar{F}^{\lambda-1}(x+t) dH(\lambda)
        \\&=&\label{dhr1} \displaystyle E\left[\frac{\bar{F}(x+t) \bar{F}^{\Lambda-1}(x+t)}{\bar{F}(x)}\right]
    \end{eqnarray}

    Now $X$ is DFR implies $\frac{\bar{F}(x+t)}{\bar{F}(x)}$ is increasing in $x$ for any $t>0$. Again  $\displaystyle\bar{F}^{\lambda -1}(x+t)$ will be increasing in $x$ for any $0<\lambda\leq 1.$ Now if we consider $\Lambda$ such that $P(0<\Lambda\leq 1)=1$ the result follows immediately.
    \item[(ii)] Similarly $X$ is IFR implies $\frac{\bar{F}(x+t)}{\bar{F}(x)}$ is decreasing in $x$. Again  $\bar{F}^{\lambda -1}(x+t)$ will be decreasing in $x$ for any $\lambda\geq 1.$ Now if we consider $\Lambda$ such that $P(\Lambda\geq 1)=1$ the result follows immediately.
\end{enumerate}
\begin{remark} Theorem \ref{th4}(i) implies that under the stated
assumptions on $X$ and $\Lambda$, $r_{X^*}(t)\leq r_X (t')$ for
$t\geq t' \geq 0$. Similarly, Theorem \ref{th4}(ii) implies that
$r_{X^*}(t)\geq r_X (t')$ for $t'\geq t \geq 0$.\end{remark}

\begin{thm}\label{th5}~
    \begin{enumerate}
        \item[(i)] $X^* \geq_{mrl\uparrow} X$ if $X$ is IMRL, provided $0<\Lambda\leq 1$ with probability 1;
        \item[(ii)] $X^* \leq_{mrl\uparrow} X$ if $X$ is DMRL, provided $\Lambda\geq 1$ with probability 1.
    \end{enumerate}
\end{thm}

\textbf{Proof:}
\begin{enumerate}
    \item [(i)] We have
 \begin{eqnarray}
\nonumber   \int_{x+t}^{\infty} \bar{F}^*(u) du / \int_{x}^\infty
\bar{F}(u)du &=&
\frac{\int_{0}^{\infty}\int_{x+t}^{\infty}\bar{F}^{\lambda}(u)dH(\lambda)}{\int_{x}^\infty
\bar{F}(u)du}
\\\label{mrl0} &=& E\left[\int_{x+t}^{\infty} \bar{F}^{\Lambda}(u) du / \int_{x}^\infty \bar{F}(u)du\right]
\end{eqnarray}
Now if $X$ is IMRL then $   \int_{x+t}^{\infty} \bar{F}(u) du /
\int_{x}^\infty \bar{F}(u)du$ increasing in $x$ for any $t>0.$ That
is we have \begin{equation}\label{mrl1}
    \frac{\bar{F}(x+t)}{\int_{x+t}^\infty \bar{F}(u) du} \leq   \frac{\bar{F}(x)}{\int_{x}^\infty \bar{F}(u) du}
\end{equation}
Let us define a function $\alpha(\lambda)=
\frac{\bar{F}^\lambda(x+t)}{\int_{x+t}^\infty \bar{F}^\lambda(u)
du}.~~~\lambda >0.$
\begin{eqnarray}
    \nonumber\alpha^\prime(\lambda) &\stackrel{sgn}{=}& \int_{x+t}^{\infty} \bar{F}^\lambda(u) [\log(\bar{F}(x+t))- \log(\bar{F}(u))] du
    \\\label{mrl2} &\stackrel{sgn}{=}& \geq 0.
\end{eqnarray}

Therefore from (\ref{mrl1}) and (\ref{mrl2}) we have for any
$0<\lambda \leq 1$ we have \begin{equation}\label{mrl3}
\frac{\bar{F}^\lambda (x+t)}{\int_{x+t}^\infty \bar{F}^\lambda (u)
du} \leq  \frac{\bar{F}(x+t)}{\int_{x+t}^\infty \bar{F}(u) du} \leq
\frac{\bar{F}(x)}{\int_{x}^\infty \bar{F}(u) du}.
\end{equation}
Hence from (\ref{mrl3}) we can easily conclude that (\ref{mrl0}) is
increasing in $x$ if $P(0 < \Lambda \leq 1 )=1 $.
\item[(ii)] Since  $X$ is DMRL hence $  \int_{x+t}^{\infty} \bar{F}(u) du / \int_{x}^\infty \bar{F}(u)du$ decreasing in $x.$ That is we have \begin{equation}\label{mrl4}
    \frac{\bar{F}(x+t)}{\int_{x+t}^\infty \bar{F}(u) du} \geq   \frac{\bar{F}(x)}{\int_{x}^\infty \bar{F}(u) du}
\end{equation}

Therefore from (\ref{mrl4}) and (\ref{mrl2}) we have for any
$\lambda \geq 1$ we have \begin{equation}\label{mrl5}
\frac{\bar{F}^\lambda (x+t)}{\int_{x+t}^\infty \bar{F}^\lambda (u)
du} \geq  \frac{\bar{F}(x+t)}{\int_{x+t}^\infty \bar{F}(u) du} \geq
\frac{\bar{F}(x)}{\int_{x}^\infty \bar{F}(u) du}.
\end{equation}
Hence from (\ref{mrl4}) we can easily conclude that (\ref{mrl0}) is
decreasing in $x$ if $ P(\Lambda \geq 1 )=1 $.
\end{enumerate}
\begin{remark} Theorem \ref{th5}(i) implies that under the stated
assumptions on $X$ and $\Lambda$, $m_{X^*}(t)\geq m_X(t')$ for
$t'\geq t \geq 0$. Similarly, Theorem \ref{th5}(ii) implies that
$m_{X^*}(t)\leq m_X (t')$ for $t\geq t' \geq 0$.\end{remark}

\begin{thm}\label{th6}~
    \begin{enumerate}
        \item[(i)] $X^* \geq_{mrl\downarrow} X$ if $X$ is DMRL, provided $0< \Lambda\leq 1$ with probability 1.
        \item[(ii)] $X^* \leq_{mrl\downarrow} X$ if $X$ is IMRL, provided $\Lambda\geq 1$ with probability 1;
    \end{enumerate}
\end{thm}

 \textbf{Proof:}
\begin{enumerate}
\item[(i)] We have
    \begin{eqnarray}
        \nonumber   \int_{x}^{\infty} \bar{F}^*(u) du / \int_{x+t}^\infty \bar{F}(u)du &=& \frac{\int_{0}^{\infty}\int_{x}^{\infty}\bar{F}^{\lambda}(u)dH(\lambda)}{\int_{x+t}^\infty \bar{F}(u)du}
        \\\label{mrl10} &=& E\left[\int_{x}^{\infty} \bar{F}^{\Lambda}(u) du / \int_{x+t}^\infty \bar{F}(u)du\right]
    \end{eqnarray}
    Now, if $X$ is DMRL then $  \int_{x+t}^{\infty} \bar{F}(u) du / \int_{x}^\infty \bar{F}(u)du$ decreasing in $x.$ That is we have \begin{equation}\label{mrl14}
        \frac{\bar{F}(x+t)}{\int_{x+t}^\infty \bar{F}(u) du} \geq   \frac{\bar{F}(x)}{\int_{x}^\infty \bar{F}(u) du}
    \end{equation}
    Therefore from (\ref{mrl14}) and (\ref{mrl12}) we have for any $0<\lambda \leq 1$ we have \begin{equation}\label{mrl15}  \frac{\bar{F}(x+t)}{\int_{x+t}^\infty \bar{F}(u) du} \geq   \frac{\bar{F}(x)}{\int_{x}^\infty \bar{F}(u) du} \geq   \frac{\bar{F}^\lambda(x)}{\int_{x}^\infty \bar{F}^\lambda(u) du}.
    \end{equation}
    Hence from (\ref{mrl14}) we can easily conclude that (\ref{mrl10}) is increasing in $x$ if $ P(0< \Lambda \leq 1 ). $
    \item [(ii)] If $X$ is IMRL then $   \int_{x+t}^{\infty} \bar{F}(u) du /
\int_{x}^\infty \bar{F}(u)du$ increasing in $x.$ That is we have
\begin{equation}\label{mrl11}
        \frac{\bar{F}(x+t)}{\int_{x+t}^\infty \bar{F}(u) du} \leq   \frac{\bar{F}(x)}{\int_{x}^\infty \bar{F}(u) du}
    \end{equation}
    Let us define a function $\beta(\lambda)=    \frac{\bar{F}^\lambda(x)}{\int_{x}^\infty \bar{F}^\lambda(u) du}.~~~\lambda >0.$
    \begin{eqnarray}
        \nonumber\beta^\prime(\lambda) &\stackrel{sgn}{=}& \int_{x}^{\infty} \bar{F}^\lambda(u) [\log(\bar{F}(x))- \log(\bar{F}(u))] du
        \\\label{mrl12} &\stackrel{sgn}{=}& \geq 0.
    \end{eqnarray}
    Therefore from (\ref{mrl11}) and (\ref{mrl12}) we have for any $\lambda \geq 1$ we have \begin{equation}\label{mrl13}   \frac{\bar{F}(x+t)}{\int_{x+t}^\infty \bar{F}(u) du} \leq    \frac{\bar{F}(x)}{\int_{x}^\infty \bar{F}(u) du}  \leq  \frac{\bar{F}^{\lambda}(x)}{\int_{x}^\infty \bar{F}^\lambda(u) du}.
    \end{equation}
    Hence from (\ref{mrl13}) we can easily conclude that (\ref{mrl10}) is decreasing in $x$ if $P( \Lambda \geq 1)=1.  $
\end{enumerate}
\begin{remark} Theorem \ref{th6}(i) implies that under the stated
assumptions on $X$ and $\Lambda$, $m_{X^*}(t)\geq m_X (t')$ for
$t\geq t' \geq 0$. Similarly, Theorem \ref{th6}(ii) implies that
$m_{X^*}(t)\leq m_X(t')$ for $t'\geq t \geq 0$.\end{remark}

\section{Results for resilience model:}
Here we study some shifted stochastic ordering of resilience models
based on some ageing properties of concerned baseline r.v.'s. Let
$X^{*}$ follow resilience model with baseline distribution $G$, and
resilience r.v. $\Omega$ having distribution function $K$ so that
the distribution function of $X^{*}$ is given by
\begin{equation}\label{rsln}
    G^*(x)= \int_{0}^{\infty}G^{\omega}(x)  dK(\omega)
\end{equation}
Throughout this section, we consider $X$ be a r.v. with distribution
function $G$ and $X^*$ be the r.v.'s as defined above for which the
distribution function is given by equation (\ref{rsln}). Also
consider that $X$ be an absolutely continuous non-negative r.v.
\par Following theorem shows that, for a baseline r.v.
$X$ with ILR (resp. DLR) property, effect of a resilience r.v.
$\Omega$ with $P(\Omega\geq 1)=1$ (resp. $P(0<\Omega\leq 1)=1$) on
$X$ is that, $X^*$ will be greater than (resp. less than) $X$ in the
sense of the up shifted likelihood ratio order.
\begin{thm}\label{th4.1}~
    \begin{enumerate}
        \item[(i)] $X^* \geq_{lr\uparrow} X$ if $X$ is ILR, provided $\Omega\geq 1$ with probability 1;
        \item[(ii)] $X^* \leq_{lr\uparrow} X$ if $X$ is DLR, provided $0<\Omega\leq 1 $ with probability 1.
    \end{enumerate}
\end{thm}
\textbf{Proof:} \begin{enumerate}

    \item[(i)] We have \begin{eqnarray}
        \nonumber \frac{g^*(x)}{g(x+t)} &=& \frac{g(x)}{g(x+t)}\times \int_{0}^{\infty} \omega G^{\omega-1}(x) dK(\omega)
        \\&=&\label{ulr5} \displaystyle E\left[\frac{g(x)\Omega G^{\Omega-1}(x)}{g(x+t)}\right]
    \end{eqnarray}

    Now $X$ is ILR implies $\frac{g(x)}{g(x+t)}$ is increasing in $x$ for any $t>0$. Again  $\omega G^{\omega -1}(x)$ will be increasing in $x$ for any $\omega\geq 1.$ Now if we consider $\Omega$ such that $P(\Omega\geq 1)=1$ the result follows immediately.
    \item[(ii)] Similarly $X$ is DLR implies $\frac{g(x)}{g(x+t)}$ is decreasing in $x$ for any $t>0$. Again  $\omega G^{\omega -1}(x)$ will be decreasing in $x$ for any $0<\omega\leq 1.$ Now if we consider $\Omega$ such that $P(0<\Omega\leq 1)=1$ the result follows immediately.
\end{enumerate}

The following corollary follows immediately in case $\Omega$ is a
degenerate r.v.
\begin{coro}~
    \begin{enumerate}
        \item[(i)] $X^* \geq_{lr\uparrow} X$ if $X$ is ILR, provided $\omega\geq 1$ ;
        \item[(ii)] $X^* \leq_{lr\uparrow} X$ if $X$ is DLR, provided $0<\omega\leq 1$ .
    \end{enumerate}
\end{coro}

\begin{thm}~
    \begin{enumerate}
        \item[(i)] $X^* \geq_{lr\downarrow} X$ if $X$ is DLR, provided $\Omega\geq 1$ with probability 1;
        \item[(ii)] $X^* \leq_{lr\downarrow} X$ if $X$ is ILR, provided $0<\Omega\leq 1 $ with probability 1.
    \end{enumerate}
\end{thm}
\textbf{Proof:} \begin{enumerate}

    \item[(i)] We have \begin{eqnarray}
        \nonumber \frac{g^*(x+t)}{g(x)} &=& \frac{g(x+t)}{g(x)}\times \int_{0}^{\infty} \omega G^{\omega-1}(x+t) dK(\omega)
        \\&=&\label{ulr6} \displaystyle E\left[\frac{g(x+t)\Omega G^{\Omega-1}(x+t)}{g(x)}\right]
    \end{eqnarray}

    Now $X$ is DLR implies $\frac{g(x+t)}{g(x)}$ is increasing in $x$ for any $t>0$. Again  $\omega G^{\omega -1}(x)$ will be increasing in $x$ for any $\omega\geq 1.$ Now if we consider $\Omega$ such that $P(\Omega\geq 1)=1$ the result follows immediately.
    \item[(ii)] Similarly $X$ is ILR implies $\frac{g(x+t)}{g(x)}$ is decreasing in $x$. Again  $\omega G^{\omega -1}(x)$ will be decreasing in $x$ for any $0<\omega\leq 1.$ Now if we consider $\Omega$ such that $P(0<\Omega\leq 1)=1$ the result follows immediately.
\end{enumerate}
Following theorem shows that, for a baseline r.v. $X$ with DRFR
(resp. IRFR) property, effect of a resilience r.v. $\Omega$ with
$P(\Omega\geq 1)=1$ (resp. $P(0<\Omega\leq 1)=1$) on $X$ is that,
$X^*$ will be greater than (resp. less than) $X$ in the sense of the
up shifted reversed hazard rate order.
\begin{thm}\label{th4.3}~
    \begin{enumerate}
        \item[(i)] $X^* \geq_{rh\uparrow} X$ if $X$ is DRFR, provided $\Omega\geq 1$ with probability 1;
        \item[(ii)] $X^* \leq_{rh\uparrow} X$ if $X$ is IRFR, provided $0<\Omega\leq 1 $ with probability 1.
    \end{enumerate}
\end{thm}
\textbf{Proof:} \begin{enumerate}

    \item[(i)] We have \begin{eqnarray}
        \nonumber \frac{G^*(x)}{G(x+t)} &=& \frac{G(x)}{G(x+t)}\times \int_{0}^{\infty}  G^{\omega-1}(x) dK(\omega)
        \\&=&\label{ulr7} \displaystyle E\left[\frac{G(x) G^{\Omega-1}(x)}{G(x+t)}\right]
    \end{eqnarray}

    Now $X$ is DRFR implies $\frac{G(x)}{G(x+t)}$ is increasing in $x$ for any $t>0$. Again  $ G^{\omega -1}(x)$ will be increasing in $x$ for any $\omega\geq 1.$ Now if we consider $\Omega$ such that $P(\Omega\geq 1)=1$ the result follows immediately.
    \item[(ii)] Similarly $X$ is IRFR implies $\frac{G(x)}{G(x+t)}$ is decreasing in $x$. Again  $ G^{\omega -1}(x)$ will be decreasing in $x$ for any $0<\omega\leq 1.$ Now if we consider $\Omega$ such that $P(0<\Omega\leq 1)=1$ the result follows immediately.
\end{enumerate}
\begin{remark} Theorem \ref{th4.3}(i) implies that under the stated
assumptions on $X$ and $\Omega$, $\tilde{r}_{X^{*}} (t)\geq
\tilde{r}_X (t')$ for $t' \geq t\geq 0$. Similarly, Theorem
\ref{th4.3}(ii) implies that $\tilde{r}_{X^{*}} (t)\leq \tilde{r}_X
(t')$ for $t \geq t'\geq 0$.\end{remark}

\begin{example}\label{ex4.1}
    Let $X$ follows Weibull r.v. with cdf $G(x)=1-e^{-2 x^2}$, $x\geq0$. Clearly, $X$ is DRFR. Let $\Omega$ to be uniformly distributed on
$[2,5].$ Then it is easy to check that $G^*(x)/G(x+t)$ is increasing
in $x$ for all $t>0.$
\end{example}

\begin{thm}\label{thmrh}~
    \begin{enumerate}
        \item[(i)] $X^* \geq_{rh\downarrow} X$ if $X$ is IRFR, provided $\Omega\geq 1$ with probability 1;
        \item[(ii)] $X^* \leq_{rh\downarrow} X$ if $X$ is DRFR, provided $0<\Omega\leq 1$ with probability 1.
    \end{enumerate}
\end{thm}
\textbf{Proof:} \begin{enumerate}

    \item[(i)] We have \begin{eqnarray}
        \nonumber \frac{G^*(x+t)}{G(x)} &=& \frac{G(x+t)}{G(x)}\times \int_{0}^{\infty}  G^{\omega-1}(x+t) dK(\omega)
        \\&=&\label{ulr8} \displaystyle E\left[\frac{G(x+t) G^{\Omega-1}(x+t)}{G(x)}\right]
    \end{eqnarray}

    Now $X$ is IRFR implies $\frac{G(x+t)}{G(x)}$ is increasing in $x$ for any $t>0$. Again  $ G^{\omega -1}(x)$ will be increasing in $x$ for any $\omega\geq 1.$ Now if we consider $\Omega$ such that $P(\Omega\geq 1)=1$ the result follows immediately.
    \item[(ii)] Similarly $X$ is DRFR implies $\frac{G(x+t)}{G(x)}$ is decreasing in $x$. Again  $\omega G^{\omega -1}(x)$ will be decreasing in $x$ for any $0<\omega\leq 1.$ Now if we consider $\Omega$ such that $P(0<\Omega\leq 1)=1$ the result follows immediately.
\end{enumerate}
\begin{remark} Theorem \ref{thmrh}(i) implies that under the stated
assumptions on $X$ and $\Omega$, $\tilde{r}_{X^{*}} (t)\geq
\tilde{r}_X (t')$ for $t \geq t'\geq 0$. Similarly, Theorem
\ref{thmrh}(ii) implies that $\tilde{r}_{X^{*}} (t)\leq \tilde{r}_X
(t')$ for $t' \geq t\geq 0$.\end{remark}

\begin{example}\label{ex4.2}
    Let $X$ follows Weibull distribution with cdf $G(x)=1-e^{-x^3}$, $x\geq0$ so that $X$ is DRFR. Let $\Omega$ to be uniformly distributed on
$[0,1].$ Then it is easy to check that $G^*(x+t)/G(x)$ is decreasing
in $x$ for all $t>0.$
\end{example}

\begin{thm}\label{thmmit}~
    \begin{enumerate}
        \item[(i)] $X^* \leq_{mit\uparrow} X$ if $X$ is IMIT, provided $\Omega\geq 1$ with probability 1;
        \item[(ii)] $X^* \geq_{mit\downarrow} X$ if $X$ is IMIT, provided $0<\Omega\geq 1$ with probability 1.
    \end{enumerate}
\end{thm}

\textbf{Proof:}
 \begin{enumerate}
\item[(i)] We have \begin{eqnarray}\label{mit0}
\nonumber\int_{0}^{x+t} G^*(u) du \Biggm/ \int_{0}^x G(u) du &=&
\displaystyle\frac{\int_{0}^{x+t}\int_{0}^\infty G^\omega(u)
dK(\omega)du }{\int_{0}^x G(u) du}
\\ &=& E \left[ \frac{\int_{0}^{x+t} G^{\Omega}(u) du}{\int_{0}^x G(u) du}\right]
\end{eqnarray}
Now $X$ is IMIT implies $\frac{\int_{0}^{x+t}G(u) du}{\int_{0}^x
G(u) du}$ decreasing in $x$ for any $t>0.$ That is we have for any
$t>0$
\begin{equation}\label{mit1}
    \frac{G(x+t)}{\int_{0}^{x+t} G(u) du} \leq  \frac{G(x)}{\int_{0}^{x} G(u) du}
 \end{equation}
Also it is easy to verify that for any $\omega >0,$
$\frac{G^{\omega}(x)}{\int_{0}^{x} G^{\omega}(u)du} $ is increasing
function of $\omega.$ Hence we have from (\ref{mit1}) for any
$0<\omega \leq 1$
\begin{equation}\label{mit3}
\frac{G^{\omega}(x+t)}{\int_{0}^{x+t} G^{\omega}(u) du}\leq
\frac{G(x+t)}{\int_{0}^{x+t} G(u) du} \leq  \frac{G(x)}{\int_{0}^{x}
G(u) du}
\end{equation}
Hence from (\ref{mit3}) we can conclude that (\ref{mit0}) is
decreasing in $x.$
\item[(ii).] Again we have
\begin{eqnarray}\label{mit5}
    \nonumber \int_{0}^{x} G(u) du \Biggm/ \int_{0}^{x+t} G^*(u) du &=& \displaystyle\frac{\int_{0}^{x} G(u) du }{\int_{0}^{x+t}\int_{0}^\infty G^\omega (u) dH(\omega) du}
    \\ &=& E \left[ \frac{\int_{0}^{x} G(u) du}{\int_{0}^{x+t} G^{\Omega}(u) du}\right]
\end{eqnarray}
As $\frac{G^{\omega}(x)}{\int_{0}^{x} G^{\omega}(u)du} $ is
increasing function of $\omega$ we have for any $\omega \geq 1$
\begin{equation}\label{mit6}
  \frac{G(x)}{\int_{0}^{x} G(u) du} \leq \frac{G(x+t)}{\int_{0}^{x+t} G(u) du}  \leq \frac{G^\omega (x+t)}{\int_{0}^{x+t} G^\omega (u) du}
\end{equation}
Consequently from (\ref{mit6}) we can conclude that (\ref{mit5}) is
decreasing in $x.$
\end{enumerate}
\begin{remark} Theorem \ref{thmmit}(i) implies that under the stated
assumptions on $X$ and $\Omega$, $mit_{X^{*}} (t)\geq mit_X (t')$
for $t \geq t'\geq 0$. Similarly, Theorem \ref{thmmit}(ii) implies
that $mit_{X^{*}} (t)\leq mit_X (t')$ for $t \geq t'\geq
0$.\end{remark}

\section{Illustration with real-world data}
Here we illustrate some of our results in two real scenarios
considering two data sets, namely ``Survival times in leukaemia" and
``Fatigue-life failures" data (\citealp{hand}). In scenario I, we
establish our results (Theorems \ref{th2} and \ref{th4}) for frailty
model. In scenario II, we establish our results for resilience model
(Theorems \ref{th4.1} and \ref{th4.3}).
\par
\begin{figure}[h]
    \begin{center}
    \includegraphics[width=10cm, height=6cm]{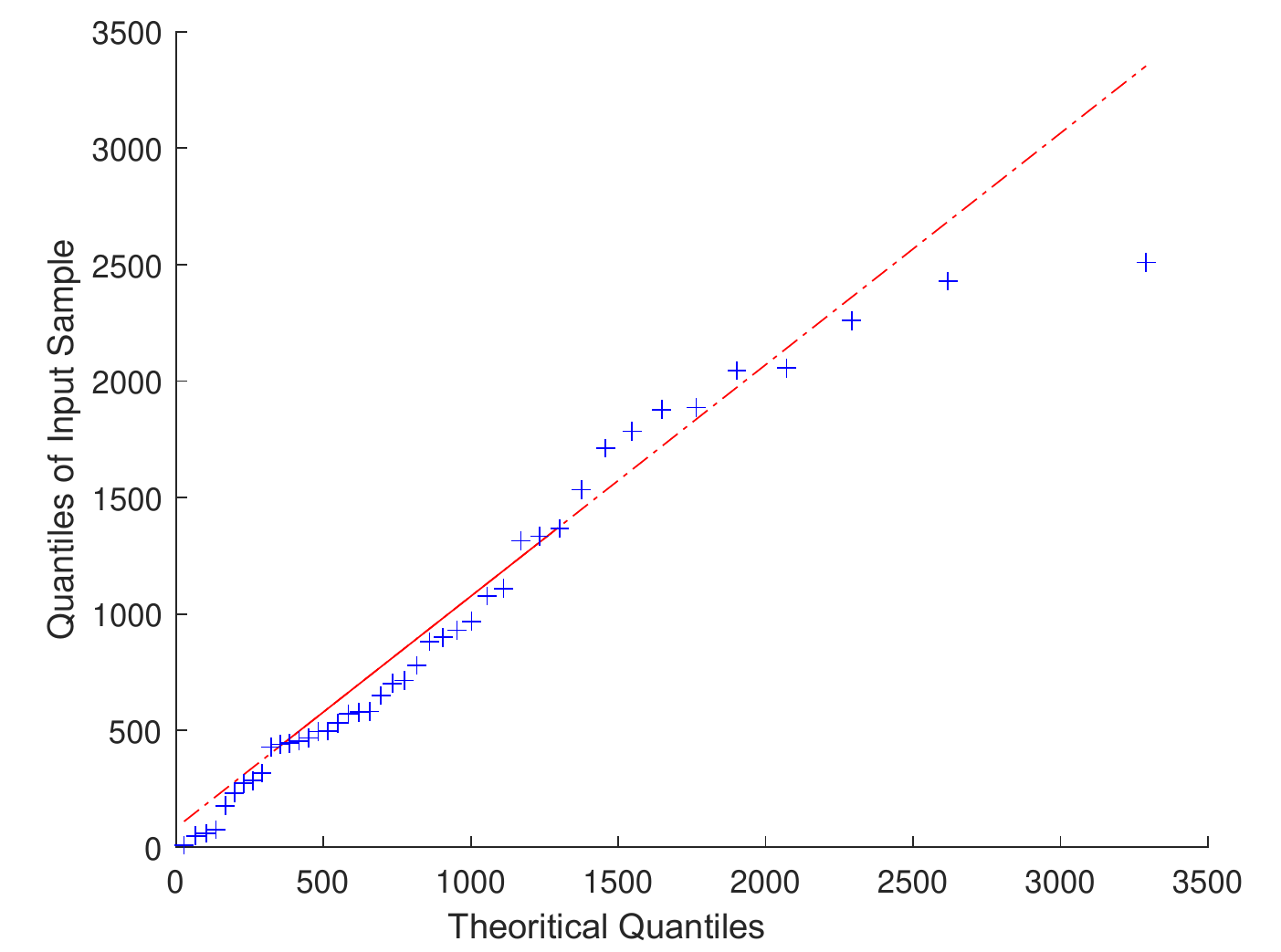}
    \caption{QQ plot of sample data vs Weibull distribution}
    \label{QQ1}
    \end{center}
\end{figure}

 \textbf{Scenario I:} We consider the data set ``Survival times in leukaemia" (\citealp{hand}) which contains the survival times of 43 patients suffering from chronic granulocytic leukaemia, measured in days from the time of diagnosis. From the quantile-quantile (Q-Q) plot (Figure \ref{QQ1}) and results of
Anderson-Darling test (Table \ref{tablead}) for the observed
samples, it is observed that Weibull distribution fits well.
Estimated values of parameters of the fitted baseline Weibull ($X$)
with the distribution function $F(x)=1-e^{(-x/\beta)^k}$, $x\geq0,
\beta>0, k>0$ are presented in Table \ref{tablewe}.
\begin{table}[h]
    \centering
    \caption{\textbf{Results of Anderson-Darling test}}
    \label{tablead}
\begin{tabular}{ ccc }
    \hline
        AD-value & p-value & Critical value(cv)\\
        \hline
    0.3616 & 0.8852 & 2.4978\\
    \hline
\end{tabular}
\end{table}
\begin{table}[h]
    \centering
    \caption{\textbf{Estimated parameters of Weibull distribution}}
    \label{tablewe}
    \begin{tabular}{ ccc }
        \hline
        Parameters & Estimated value & $95\%$ confidence interval\\
        \hline
        Scale $(\beta)$ & $986.672$  & $ [766.52, 1270.06]$ \\
        Shape $(k)$ & $1.24044$    &  $[0.973535, 1.58052]$
        \\
        \hline
    \end{tabular}
\end{table}

With shape parameter $k>1$, this baseline Weibull distribution is
ILR and and so is IFR. Next we consider well known Gamma-frailty
i.e. $\Lambda \sim \Gamma(1/a^2,1/a^2)$ where $\Lambda \geq 1$ with
probability 1. According to Theorem \ref{th2}(ii), the effect of
considered gamma frailty on $X$ is that, $X^* \leq_{lr\downarrow}
X$, which implies that $\kappa_{X^*}(t)\leq \kappa_X (t')$ for
$t'\geq t \geq 0$. Similarly, according to Theorem \ref{th4}(ii),
$X^* \leq_{hr\downarrow} X$, which implies that $r_X^{*} (t)\geq r_X
(t')$ for $t' \geq t\geq 0$. Also, we have $X^* \geq_{disp} X$,
where `disp' stands for dispersive order (\citealp{shaked}). It
follows from the fact that for two non-negative r.v.'s $X$ and $Y$,
$X \leq_{hr\downarrow} Y \Rightarrow X \leq_{disp} Y$
(\citealp{lillo}).
\par
\begin{figure}[h]
    \begin{center}
    \includegraphics[width=10cm, height=7cm]{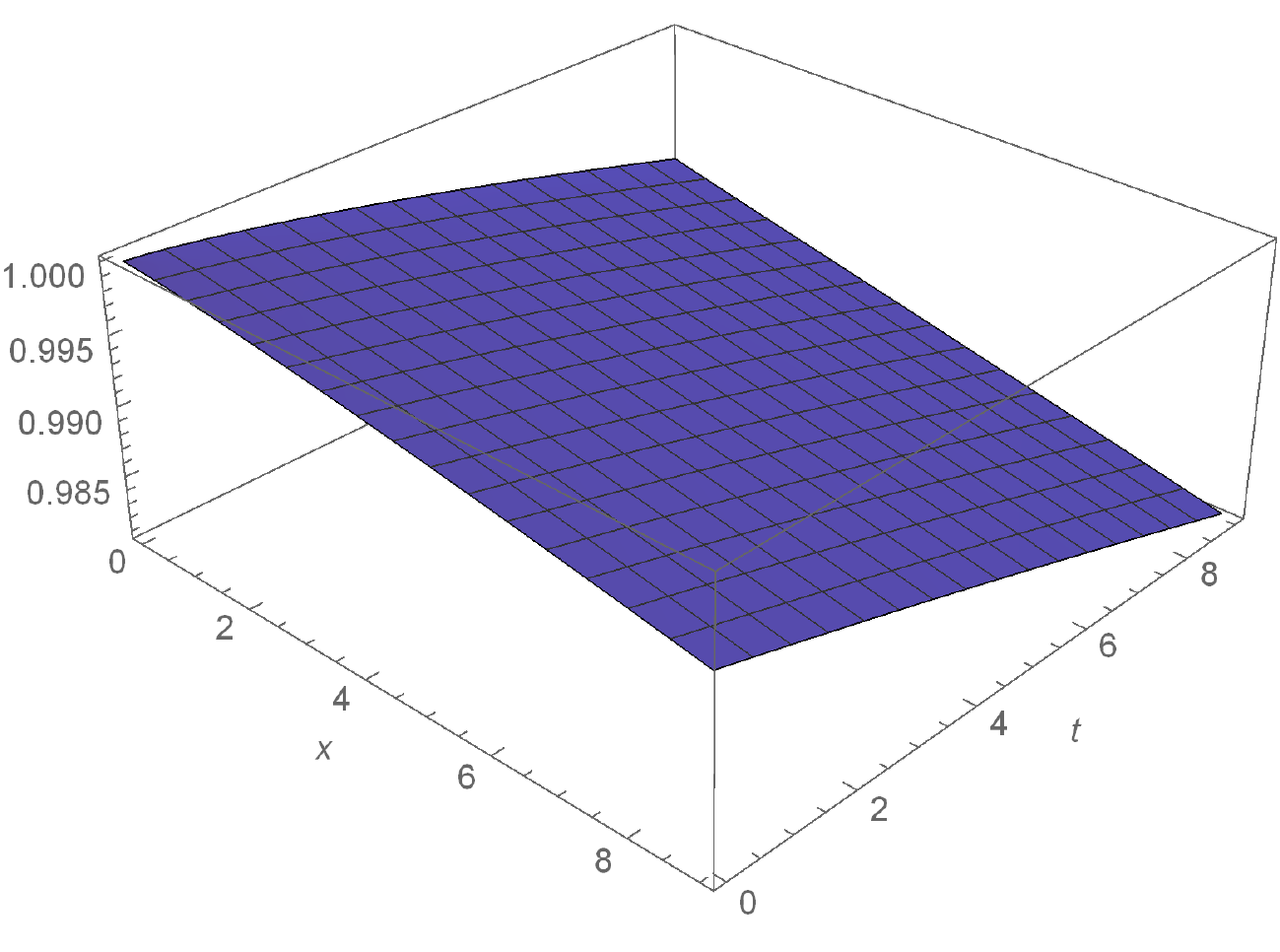}
    \caption{Plot of $f^*(x+t)\Big/f(x)$}
    \label{fltlr}
    \end{center}
\end{figure}

\begin{figure}[h]
    \begin{center}
    \includegraphics[width=10cm, height=6cm]{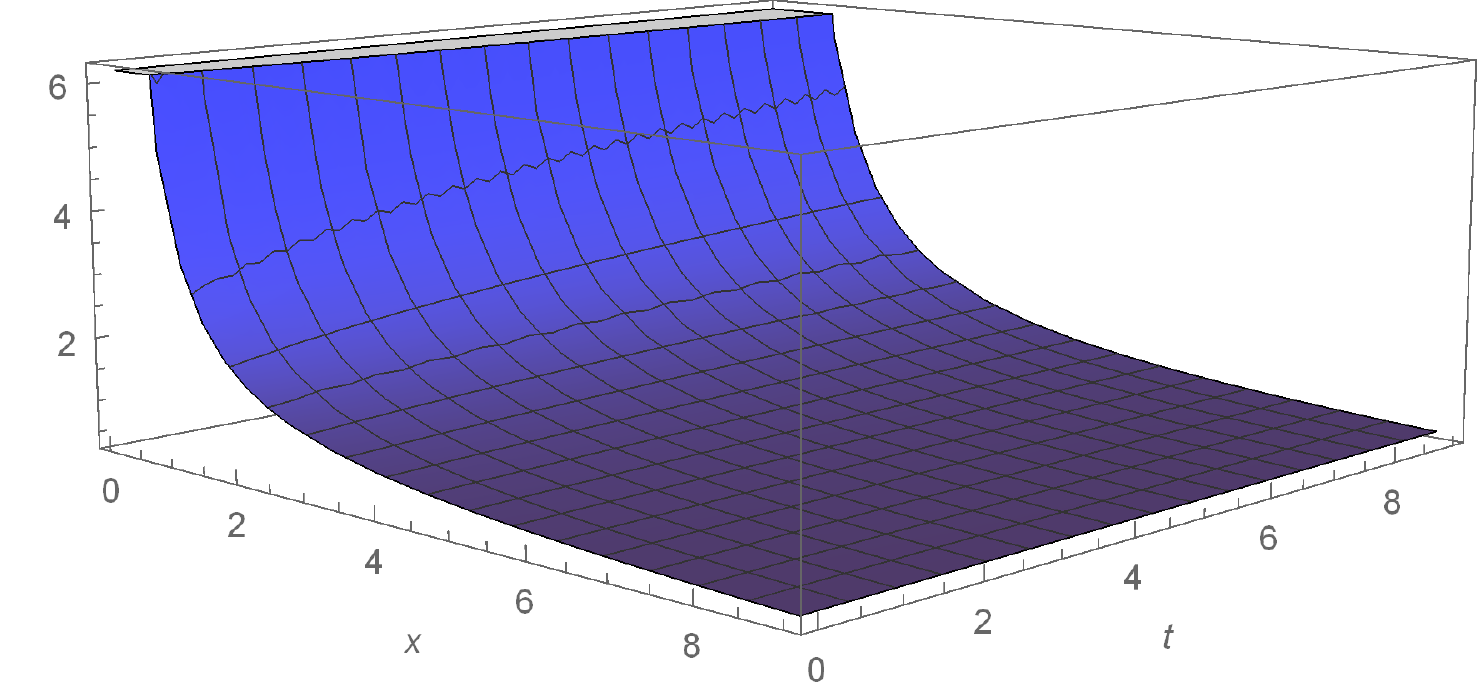}
    \caption{Plot of $\bar{F}^*(x+t)\Big/\bar{F}(x)$}
    \label{flthr}
    \end{center}
\end{figure}

 To demonstrate the above mentioned stochastic orders, we
proceed as follows. The survival function and probability density
function of the above frailty model (Gamma-frailty Weibull-baseline)
are, respectively
\begin{eqnarray}
    \bar{F}^*(t) &=&  \frac{(1/a^2)^{-1/a^2}(a^2 + \frac{t^k}{\beta^k})^{-1/a^2} \zeta_1\left(\frac{1}{a^2},a^2+\frac{t^k}{\beta^k}\right)}{\Gamma(\frac{1}{a^2}) \left(  1- \zeta_2 (\frac{1}{a^2},0,a^2) \right)}
    \\ \text{and}~ f^*(t) &=&    \frac{\frac{k t^{k -1} }{\beta^k }   (1/a^2)^{-1/a^2}  (a^2 + \frac{t^k}{\beta^k})^{-1-1/a^2} \zeta_1\left(\frac{1}{a^2},a^2+\frac{t^k}{\beta^k}\right)}{\Gamma(\frac{1}{a^2}) \left(  1- \zeta_2 (\frac{1}{a^2},0,a^2) \right)}\end{eqnarray}
Let $t_1,t_2,...,t_n$ be the observations under consideration. We
now obtain maximum likelihood estimation of the parameter $a$ under
the Gamma-frailty Weibull-baseline. The likelihood function is given
by
    \begin{eqnarray} \mathcal{L}(a| t_1,t_2,..,t_n) &=& \left(\frac{(\frac{1}{a^2})^{-1/a^2}}{\Gamma(\frac{1}{a^2}) \left(  1- \zeta_2 (\frac{1}{a^2},0,a^2) \right)} \right)^n k^n \prod_{i=1}^n \left( t_i \right)^{k-1} \prod_{i=1}^n (a^2 + \frac{t_i^k}{\beta^k})^{-1-1/a^2}
    \\&& ~~~~~~~~ \times \prod_{i=1}^n
    \zeta_1\left(\frac{1}{a^2},a^2+\frac{t_i^k}{\beta^k}\right),
\end{eqnarray}
where $\zeta_1(a,x) = \int_{x}^\infty t^{a-1} e^{-t} dt$ and
$\zeta_2(a,x) = \frac{\int_{0}^x t^{a-1} e^{-t} dt}{\Gamma(a)}$ are
upper incomplete gamma functions and regularized lower incomplete
gamma functions respectively. Estimated value of $a$ is obtained as
$0.784$ with
$\mathbb{P}(\Gamma \geq 1) = 1$.\\
We then plotted $f^*(x+t)\Big/f(x)$ taking some finite range of $x$
and $t$ as shown in Figure \ref{fltlr}, which is clearly showing
that the ratio is decreasing in $x$, giving $X^*\leq_{lr\downarrow}
X$. To demonstrate that $X^* \leq_{hr\downarrow} X$, we plotted
$\bar{F}^*(x+t)\Big/\bar{F}(x)$ in Figure \ref{flthr} showing that
it is decreasing in $x$.

\par
\begin{figure}[h]
    \begin{center}
        \includegraphics[width=10cm, height=6cm]{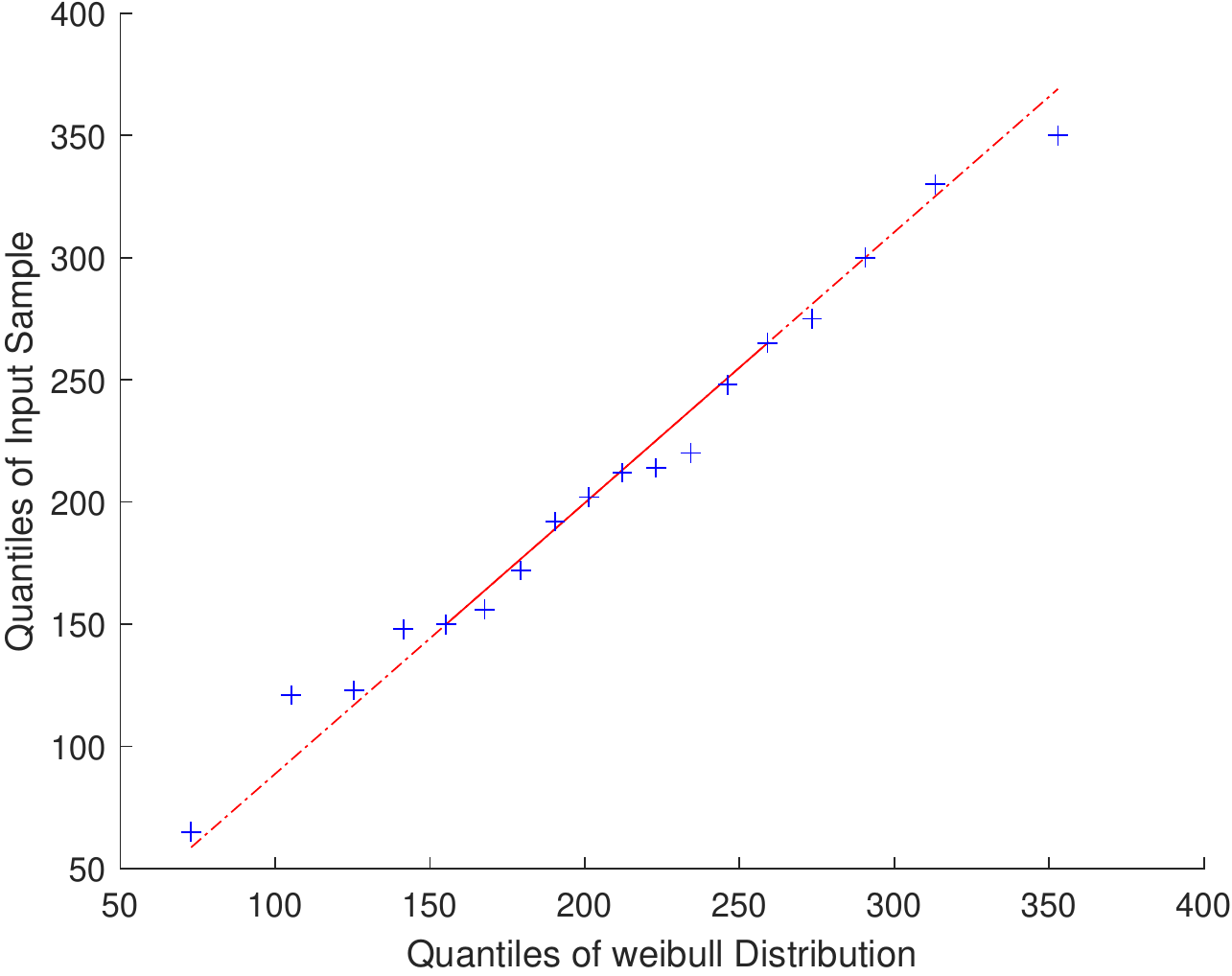}
        \caption{QQ plot of sample data vs Weibull distribution}
        \label{QQR}
    \end{center}
\end{figure}

 \textbf{Scenario II:} Here we consider the data set ``Fatigue-life failures" (\citealp{hand}) on the fatigue-life failures of ball-bearings.
The data give the number of cycles to failure. From the
quantile-quantile (Q-Q) plot (Figure \ref{QQR}) and the results of
Anderson-Darling test (Table \ref{tablead2}) for the observed
samples, it is observed that the samples can taken to be from
Weibull distribution. Estimated values of parameters of baseline
Weibull with the distribution function $G(x)=1-e^{(-x/\beta)^k}$,
$x\geq0, \beta>0, k>0$ are given in Table \ref{tablewe2}.
\begin{table}[h]
    \centering
    \caption{\textbf{Results of Anderson-Darling test}}
    \label{tablead2}
\begin{tabular}{ ccc }
    \hline
        AD-value & p-value & Critical value(cv)\\
        \hline
    0.1496 & 0.99 & 2.503\\
    \hline
\end{tabular}
\end{table}
\begin{table}[h]
    \centering
    \caption{\textbf{Estimated parameters of Weibull distribution}}
    \label{tablewe2}
    \begin{tabular}{ ccc }
        \hline
        Parameters & Estimated value & $95\%$ confidence interval\\
        \hline
        Scale & $232.9$  & $[198.758, 272.906] $ \\
        Shape & $3.0721$    &  $[2.13732, 4.41572]$
        \\
        \hline
    \end{tabular}
\end{table}

\begin{figure}[h]
    \begin{center}
    \includegraphics[width=10cm, height=7cm]{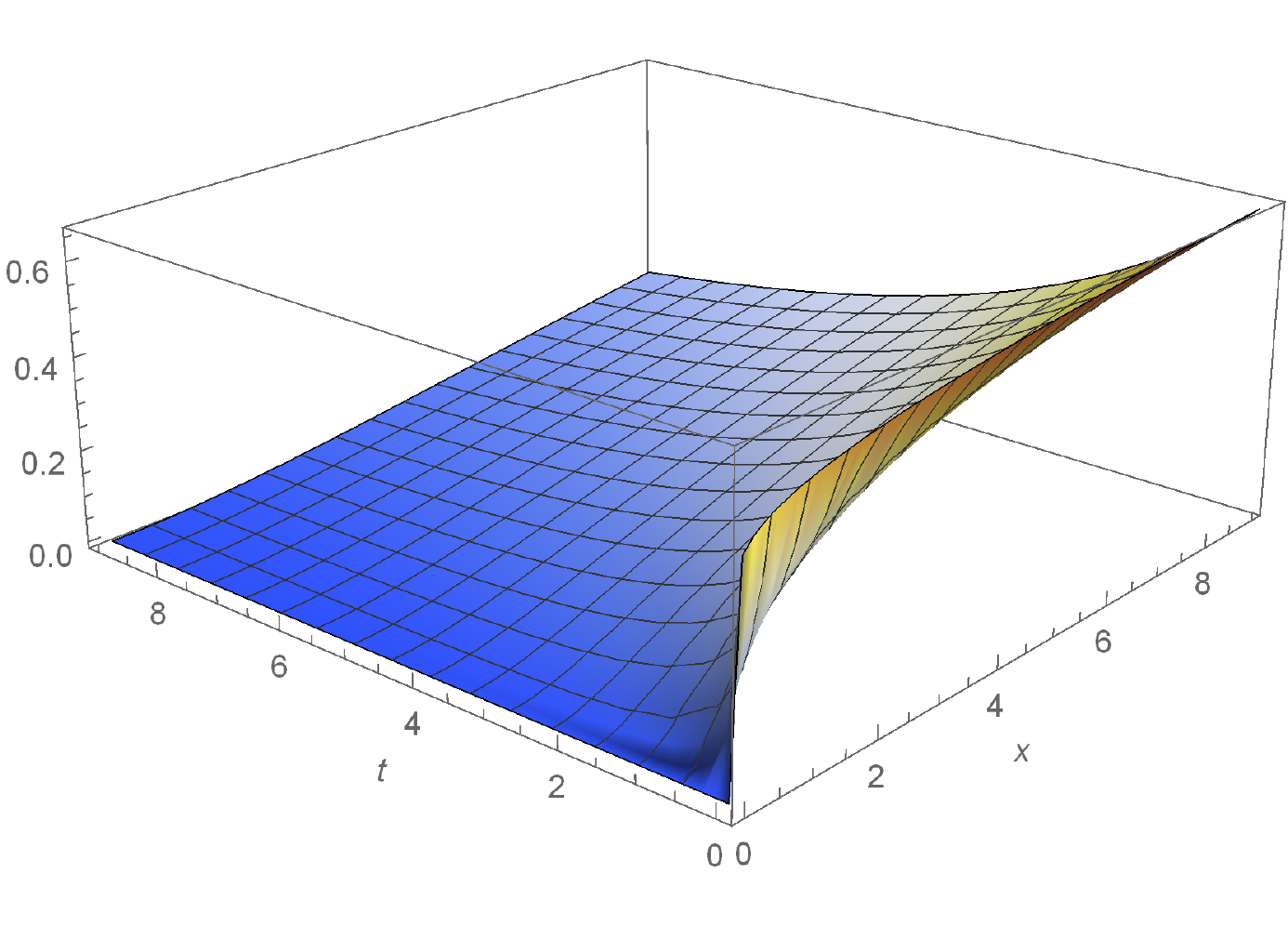}
    \caption{Plot of $g^*(x)\Big/g(x+t)$}
    \label{rslru}
    \end{center}
\end{figure}

With shape parameter $k>1$, this baseline Weibull distribution is
ILR and also is DRFR. Next we consider Gamma resilience i.e. $\Omega
\sim \Gamma(1/a^2,1/a^2)$ where $\Omega\geq 1$ with probability 1.
According to Theorem \ref{th4.1}(i), the effect of considered gamma
resilience on $X$ is that, $X^* \leq_{lr\downarrow} X$. Similarly,
according to Theorem \ref{th4.3}(i), $X^*\leq_{rh\uparrow} X$, which
indicates that $\tilde{r}_X^{*} (t)\leq \tilde{r}_X (t')$ for $t
\geq t'\geq 0$.
\par To demonstrate the above mentioned stochastic orders, we
proceed as follows. The distribution function and probability
density function of the above resilience model (Gamma-resilience
Weibull-baseline) are, respectively
\begin{eqnarray*}
    G^*(t) &=&  \frac{(1/a^2)^{-1/a^2}(a^2 - \ln(1-e^{(\frac{t}{\beta})^k}))^{-1/a^2} \zeta_1\left(\frac{1}{a^2},a^2 - \ln(1-e^{(\frac{t}{\beta})^k})\right)}{\Gamma(\frac{1}{a^2}) \left(  1- \zeta_2 (\frac{1}{a^2},0,a^2) \right)}
    \\ \text{and}~ g^*(t) &=&    \frac{\frac{k t^{k -1} }{\beta^k }   (1/a^2)^{-1/a^2}  (a^2 - \ln(1-e^{(\frac{t}{\beta})^k}))^{-1-1/a^2} \zeta_1\left(\frac{1}{a^2},a^2- \ln(1-e^{(\frac{t}{\beta})^k})\right)}{\Gamma(\frac{1}{a^2}) \left(  1- \zeta_2 (\frac{1}{a^2},0,a^2) \right)}\end{eqnarray*}

    Let $t_1,t_2,...,t_n$ be the observations under consideration. We now obtain maximum likelihood estimate of the parameter $a$ under the Gamma-resilience
Weibull-baseline. The likelihood function is given by
    \begin{eqnarray*} \mathcal{L}(a| t_1,t_2,..,t_n) &=& \left(\frac{(\frac{1}{a^2})^{-1/a^2}}{\Gamma(\frac{1}{a^2}) \left(  1- \zeta_2 (\frac{1}{a^2},0,a^2) \right)} \right)^n \prod_{i=1}^n \left( t_i \right)^{k-1} \prod_{i=1}^n (a^2 - \ln(1-e^{(\frac{t_i}{\beta})^k}))^{-1-1/a^2}
    \\&& ~~~~~~~~ \times \prod_{i=1}^n  \zeta_1\left(\frac{1}{a^2},a^2 - \ln(1-e^{(\frac{t_i}{\beta})^k})
    \right),
\end{eqnarray*}
where $\zeta_1(a,x) \text{ and } \zeta_2(a,x)$ are defined in
previous case. Estimated value of the parameter $a$ is obtained as
$4.0558$ with $\mathbb{P}(\Omega \geq 1) = 1$.
\par
\begin{figure}[H]
    \begin{center}
    \includegraphics[width=10cm, height=6cm]{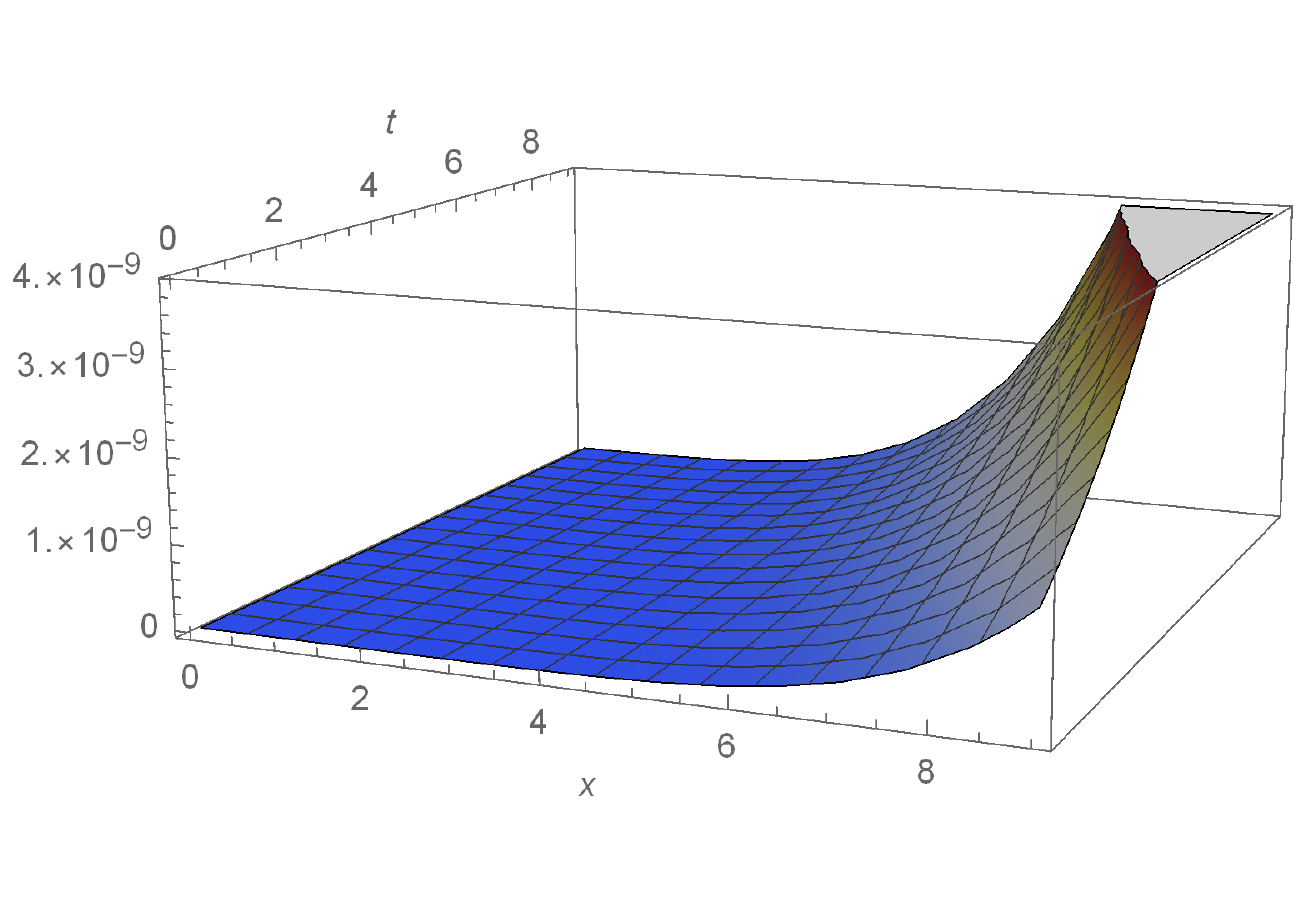}
    \caption{Plot of $G^*(x)\Big/G(x+t)$}
    \label{rsrhu}
    \end{center}
\end{figure}

 Then we plotted $g^*(x)\Big/g(x+t)$ taking some finite range of $x$ and $t$ as shown
in Figure \ref{rslru}, which is clearly showing that the ratio is
increasing in $x$, giving $X^*\leq_{lr\uparrow} X$. To demonstrate
that $X^*\leq_{rh\uparrow} X$, we plotted $G^*(x)\Big/G(x+t)$ in
Figure \ref{rsrhu} showing that it is increasing in $x$.
\section{Conclusion}
In this study, we have derived results on stochastic comparisons for
frailty as well as resilience models to study the effects of frailty
and resilience r.v.'s on the baseline r.v.'s based on some ageing
properties of concerned baseline r.v.'s. To derive the results we
have used some shifted stochastic orders which are stronger than
their respective usual counterparts, and also provide more
flexibility in stochastic comparisons. As a future study,
comparisons for considered frailty or resilience models could be
explored using other generalized stochastic orders like proportional
stochastic and shifted proportional stochastic orders.

\section*{Conflict of interest}

On behalf of all authors, the corresponding author states that there
is no conflict of interest.

\section*{Data availability statement}
The datasets analysed during the current study are available in
\cite{hand}.

\end{document}